\documentclass{article}
\usepackage{epic,latexsym,amssymb}
\usepackage{color}
\usepackage{tikz}
\usepackage{amsfonts,epsf,amsmath,subcaption,caption}

\newtheorem{thm}{Theorem}
\newtheorem{lem}{Lemma}

\newtheorem{cor}{Corollary}

\newtheorem{problem}{Problem}

\let\oldenumerate\enumerate
\renewcommand{\enumerate}{
  \oldenumerate
  \setlength{\itemsep}{0pt}
  \setlength{\parskip}{0pt}
  \setlength{\parsep}{0pt}
}

\begin{document}

\title{Drunk Angel and Hiding Devil}
\markright{drunk Angel and Devil by Normal Distribution}

\author{Nuttanon Songsuwan\thanks{Research supported by Petchra Pra Jom Klao Ph.D. Research Scholarship from King Mongkut's University of Technology Thonburi~(33/2563).}, \, Anuwat Tangthanawatsakul
\\ Pawaton Kaemawichanurat
\\ \\
Mathematics and Statistics with Applications~(MaSA),\\
Department of Mathematics, Faculty of Science,\\ King Mongkut's University of Technology Thonburi, \\
Bangkok, Thailand \\
\small \tt Email: nuttanon.19701@gmail.com, anuwat.sae@kmutt.ac.th, pawaton.kae@kmutt.ac.th}

\date{}
\maketitle

\begin{abstract}
The angel game is played on $2$-dimensional infinite grid by $2$ players, the angel and the devil. In each turn, the angel of power $c \in \mathbb{N}$ moves from her current point $(x, y)$ to a point $(x', y')$ which $\max\{|x - x'|, |y - y'|\} \leq c$ while the devil chooses a point to destroy in his turn. Then, the angel can no longer land on these destroyed points. The angel wins if she has a strategy to escape from the devil forever and the devil wins if he can cage the angel in his destroyed points by a finite number of turns. It was proved in 2007 that the angel of power at least $2$ always wins. In this paper, we rise the problem when the angel is drunk. She randomly moves to any point in the range of her power in each turn. In our game version, the devil must cage the angel by a given finite number of turns, otherwise, the angel wins. We present a strategy for the devil that: if the devil plays with this strategy, then for given $c \in \mathbb{N}$ and $\epsilon > 0$, the devil can cage the angel of power $c$ with probability greater than $1 - \epsilon$ if and only if the game is played on an $n$-dimensional infinite grid when $n \leq 2$. We also establish the results related to the hitting time once the angel is first time outside an $n$-dimensional sphere of a given radius. The numerical simulation results are also presented in the last section.
\end{abstract}


\section{Introduction and Motivation}
For a positive integer $n$, a graph $G$ is said to be an $n$\emph{-dimensional infinite grid} if every vertex $x$ of $G$ is an $n$-tuple $x = (x_1,\dots , x_n)$ where $x_1, x_{2}, ..., x_{n} \in \mathbb{Z}$ and two vertices $x$ and $y$ of $G$ are adjacent if $\sum_{i=1}^{n}|x_i - y_i| = 1$. For an $n$-dimensional infinite grid graph, we denote the \textit{supremum norm} of vertex $x$ by $\|x\| = \max_{1 \le i \le n} |x_i|$. The \emph{distance} between $x$ and $y$ is defined by $d(x, y) = \|x - y\|$. The \emph{neighbor set} $N(x)$ of the vertex $x$ is the set of all vertices that are adjacent to $x$. The \emph{close neighbor set} of $x$ is the set $N(x) \cup \{x\}$ and is denoted by $N[x]$. In an $n$-dimensional infinite grid graph, the \emph{sphere} of radius $k$ centered at the origin is the set of all points whose distance from the origin is at most $k$ and is denoted by $S^{n}_{k}$. The \emph{lattice sphere} of radius $k$ centered at the origin is the set of all grid points (vertices) whose distance from the origin is at most $k$ and is denoted by $L^{n}_{k}$. Further, the \emph{hollow lattice sphere} of inner radius $k$ with the thickness $c$ and centered at the origin is $H^{n}_{k, c} = L^{n}_{k + c} - L^{n}_{k - 1}$. 

\indent A random walk is a mathematical process which is used to describe a random move of an object in a mathematical space. Over the past one hundred years, this method has been dramatically evolved and applied to many scientific areas. An elementary example of a random walk study is when it takes place on graphs. For a graph $G$, a random walk on $G$ can be defined as the Markov chain $\{X_{k}, k \geq 0\}$ where each step moves from the vertex $x$ to the vertex $y$ with the probability $\frac{1}{|N[x]|}$ if $y \in N[x]$ and with the probability $0$ if $y \notin N[x]$. Random walks on graphs have been used to present structures of computer and electrical networks. For examples, see \cite{bartal2002more,curado2022anew,doyle1984random,gatto2018saddlepoint,isler2005randomized,kang2004random,koolen2016collection,nash1959random,novikov2020random,tetali1991random,wang2018expected}.


The angel game appeared in Berlekamp et. al \cite{berlekamp2018winning} and has been popularized by Conway \cite{conway1996angel} in 1996. In the classical version, the game is played on a $2$-dimensional infinite grid graph, $\mathbb{Z}^2$. It consists of 2 players, the angel and the devil, who play alternately. Initially, the angel lands on some vertex, $(0,0)$ say. In her turn, she leaps to another vertex at distance at most $c$, the \emph{power} of the angel, from her current position. In each turn of the devil, he destroys one vertex. Afterwards, the angel cannot land on this vertex during the game. The angel loses if she is caged in a closed shape accumulated by the destroyed vertices. The angel wins if she has a strategy to escape forever. Both players play optimally. It is known that for $c = 1$, the devil wins \cite{berlekamp2018winning}. In 1996, Conway \cite{conway1996angel} proved that if the angel always moves upward (the angel moves from $(x,y)$ to $(x^*,y^*)$ for which $y^* \ge y$), then the devil wins. For $c \ge 2$, the proof has been presented independently by Kloster \cite{kloster2007solution} and M\'ath\'e \cite{mathe2007angel} that the angel of power at least 2 wins.

For other variations of angel and devil, Kutz \cite{kutz2005conway} showed that, when the game is played on a $3$-dimensional infinite grid graph, the angel wins if her power is at least $13$. Bollob\'as and Leader \cite{bollobas2006angel} independently applied probabilistic proof to show that the angel wins if her power is large enough. Kutz and P\'or \cite{kutz2005angel} generalized the game to when the power of the angel is a positive real number. 

In this paper, we consider when the angel is drunk. In other words, on each of her turns, she performs a symmetric random walk on an $n$-dimensional infinite grid graph. As the angel plays randomly, the advantage belongs to the devil because this is an infinite game. So, we may rule the devil to have only one chance to predict where the angel will be after a number of moves. The precise definition of the game will be given in the next section.

\subsection{Rule of the game}

Before the game starts, the \emph{hiding} devil selects a natural number $N$. It can be extremely large but it must be finite. It will be the number of turns that the \emph{drunk} angel can move (randomly) as well as the number of vertices that the hiding devil has to destroy to cage the drunk angel. In each turn, the angel at a vertex $x = (x_{1}, ..., x_{n}) \in \mathbb{Z}^n$ moves to a vertex $y = (y_{1}, ..., y_{n})$ which $d(x, y) \leq c$ with probability $\frac{1}{(2c + 1)^{n}}$. For the devil, he \emph{secretly} destroys one vertex in each of the Turns $1$ to $N$. During these turns, the destroyed vertices do not effect the angel-move probability. She, the angel, can even moves on these vertices. However, after the angel has moved in Turn $N$, the devil suddenly reveals a closed $n$-dimensional shape of thickness $c$ (the minimum distance between the inner and outer parts of the shape is equal to $c$) that is made of the destroyed vertices. If the angel is not in the inner part of the close shape, she wins. Otherwise, the devil wins.


\indent From the rule of the game, the problem that is arisen is:

\begin{problem}
For given $c, n \in \mathbb{N}$ and real number $\epsilon > 0$, does there exist a strategy for the hiding devil to cage the drunk angel of power $c$ in an $n$-dimensional infinite grid graph with probability greater then $1 - \epsilon$ ?
\end{problem}

\indent This paper is organized as follows. The devil strategy and main results are presented in Section \ref{sec2} while the proofs are given in Sections \ref{sec4} - \ref{sec7}. We provide some related tools in Section \ref{sec3} and some numerical results in Section \ref{numerical}.

\section{The Devil Strategy and Main results}\label{sec2}
We present our main theorems in this section while the numerical results for implementation are given in Section \ref{numerical}. The first theorem shows the difficulty when we find the probability that the drunk angel will be inside $H^{1}_{k, c}$ after a number of random moves.

\begin{lem}\label{dim}
In $1$-dimensional infinite grid graph, if the probability that the drunk angel moves from her current vertex $x$ to a vertex $y$ which $d(x, y) \leq c$ is $\frac{1}{2c + 1}$, then the probability that she is at a vertex with the distance at most $k$ from $x$ after moving $l$ turns is
\begin{equation*}
\frac{1}{(2c + 1)^{l}}\sum_{m=0}^{l} (-1)^m\binom{l}{m}\binom{l-1+k+cl-m(2c+1)}{l-1}.
\end{equation*}
\end{lem}

\noindent By considering the large expression of the result of Lemma \ref{dim}, we see that there will also be a huge computation to find the probability that the drunk angel will be inside $H^{n}_{k, c}$ when the number of dimension is increased. Hence, we may apply de Moivre–Laplace theorem (see Subsection \ref{CLT}) to find this probability. 
\vskip 5 pt

\subsection{The devil strategy}
\indent As mentioned in the game rule, the hiding devil keeps (secretly) destroying the vertices to accumulate an $n$-dimensional closed shape of thickness $c$, $B^{n}_{c}$,to cage the drunk angel of power $c$ in the inner part of $B^{n}_{c}$. The question is: 
\vskip 5 pt

\indent what does the shape $B^{n}_{c}$ look like to maximize the caging probability? 
\vskip 5 pt

\begin{figure}[h]
\centering
\subcaptionbox{$c = 1$}
{\includegraphics[width=0.33\textwidth]{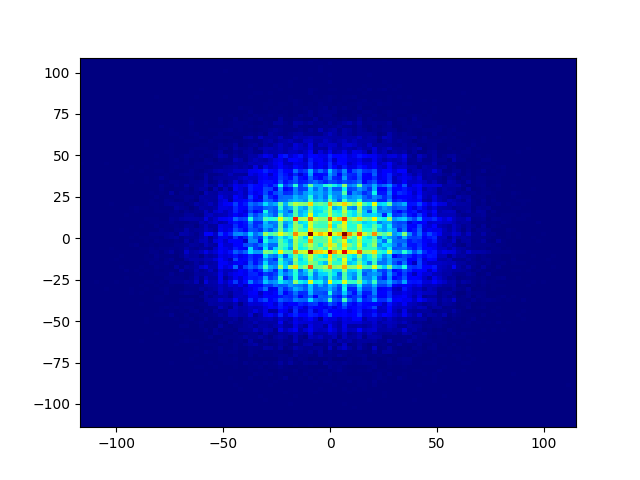}}%
\hfill
\subcaptionbox{$c = 4$}
{\includegraphics[width=0.33\textwidth]{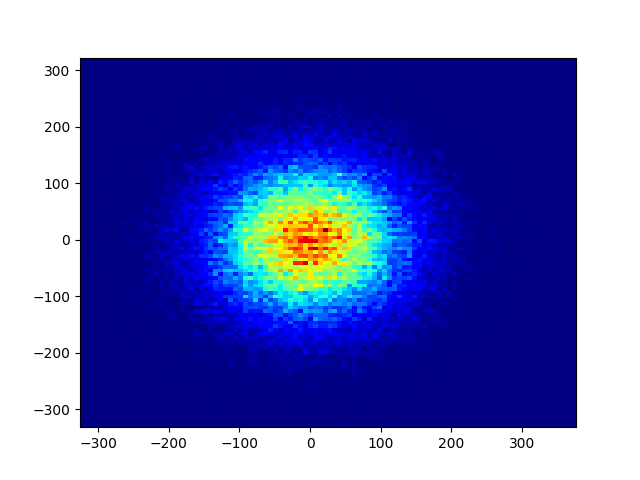}}%
\hfill
\subcaptionbox{$c = 9$}
{\includegraphics[width=0.33\textwidth]{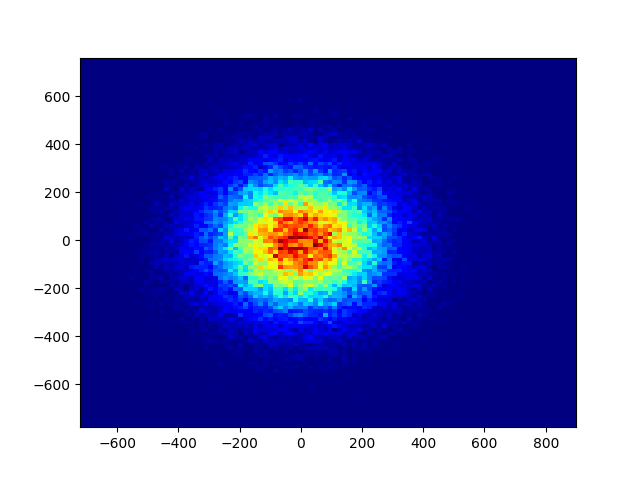}}%
\caption{The foot prints of the angels of powers $1, 4$ and $9$ after 1000 steps}
\label{fig:foot}
\end{figure}

\noindent To answer this question, we have simulated a computer programming to find what the angel footprint look like. Our results of the angel footprint can be seen in https://github.com/nuttanon19701/DAnHD, all of which illustrate a similar shape which is the circle (2D-Space). Figure \ref{fig:foot} illustrates some of our simulation results too.

\noindent More solid tool to answer this question is to employ central limit theorem. We may let $\boldsymbol{\xi}_i = (\xi_{i}^{(1)},\xi_{i}^{(2)},\dots,\xi_{i}^{(n)})$ be an $n$-dimensional random vector with $\xi_{i}^{(j)}$ be the discrete uniform random variable from $-c,-c+1,\dots,0,\dots,c-1,c$ for all $j = 1,2,\dots,n$. The position of the angel at step $m$ (starting at the origin) can be represented by the sum of random vector $\boldsymbol{\xi}_{i}$ from $i = 0,1,\dots,m$, denoted by
\begin{equation*}
    \mathbf{X}_m = \sum_{i=0}^{m} \boldsymbol{\xi}_{i}.
\end{equation*}
\noindent By multivariate central limit theorem (see Theorem \ref{MCLT}), the random vector  $\frac{\mathbf{X}_m}{\sqrt{m}}$ converges in distribution to the normal distribution $\mathcal{N}_n(\mathbf{0},\boldsymbol{\Sigma})$. This convinces that the position of the angel at the last turn is most likely symmetrically around the origin, fitting inside a hollowed sphere with appropriate inner radius.
\vskip 5 pt

\indent Hence, by this arguments and the simulations, we believe that, for any given natural numbers $k$ and $c$, the $n$-dimensional closed shape of thickness $c$ with the shortest distance from the origin to any vertex of its body is equal to $k$ that maximizes the caging probability is $H^{n}_{k, c}$. So, for the game that plays on an $n$-dimensional infinite grid graph with the drunk angel of power $c$, our game strategy for the hiding devil is to construct $H^{n}_{k, c}$ with an appropriate $k$. By the Gauss's circle problem (see Subsection \ref{gauss}), the hiding devil needs 

\begin{equation*}
N = |H^{n}_{k, c}| \approx \frac{\pi^{\frac{n}{2}}[(k+c\sqrt{n})^n]}{\Gamma(\frac{n}{2}+1)} - \frac{\pi^{\frac{n}{2}}[k^n]}{\Gamma(\frac{n}{2}+1)}
\end{equation*}

\noindent turns to complete his cage construction where $\Gamma$ is the well-known gamma function. Note that, this is the number of turns that the drunk angel can move too. We may call our strategy \emph{Gauss's circle}. 
\vskip 5 pt

\subsection{Main results}
\indent Throughout this paper, let $A^{n}_{k, c}$ be the events that the angel of power $c$ is caged by the devil after $N$ moves when the hollow sphere $H^{n}_{k, c}$is exposed. In Theorems \ref{submain1}, \ref{submain2} and \ref{submain3}, the probabilities of $A^{n}_{k, c}$ depending on $n$, the number of dimensions, are presented. The proofs of which are also given in Sections \ref{sec4} - \ref{sec6}.

\begin{thm}\label{submain1}
If the hiding devil plays with the Gauss's circle strategy on an $n$-dimensional infinite grid graph when $n \leq 2$, then, for a given $c \in \mathbb{N}$ and a real number $\epsilon > 0$, there exists $k \in \mathbb{N}$ such that $P[A^{n}_{k, c}] > 1 - \epsilon$.
\end{thm}

\begin{thm}\label{submain2}
If the hiding devil plays with the Gauss's circle strategy on a $3$-dimensional infinite grid graph, then, for a given $c \in \mathbb{N}$, we have that, for all $k \in \mathbb{N}$,
\begin{equation*}
0 < l(c) \leq P[A^{3}_{k, c}] \leq u(c)
\end{equation*}
\noindent where $l(c) = exp\left(-\frac{\sqrt{3}}{8\pi c^2(c+1)}\right)\left(\frac{1}{2c(c+1)}\right)^3\left(\frac{2\sqrt{3}}{3\pi^4(c+1)^3}\right)^{\frac{1}{2}}$ and \\
$u(c) = exp\left( -\frac{\sqrt{3}}{8\pi c^2(c+1)^3}\right)\left(\frac{1}{8\sqrt{3}\pi^2 c^2(c+1) }\right)^{\frac{1}{2}}\left(\frac{1}{2\pi c^2(c+1) - \sqrt{3}}\right)$.
\end{thm}

\begin{thm}\label{submain3}
If the hiding devil plays with the Gauss's circle strategy on an $n$-dimensional infinite grid graph when $n \geq 4$, then, for any given $c \in \mathbb{N}$, we have that
\begin{equation*}
P[A^{n}_{k, c}] \rightarrow 0 \text{ when } k \rightarrow \infty.
\end{equation*}
\end{thm}

\indent In Theorems \ref{stopping} and \ref{prob}, we apply martingale to establish the hitting time once the angel is first time outside the $n$-dimensional sphere of the given radius. The proofs of which are given in Section \ref{sec7}.

\begin{thm}\label{stopping}
Let $k$ be the positive integer. If the angel of power $c$ plays on an $n$-dimensional infinite diagonal grid graph, then
\begin{equation*}
    \frac{3}{nc(c+1)}k^2 \le  \mathbb{E}(\tau_k) \le \frac{3}{nc(c+1)}(k+c\sqrt{n})^2
\end{equation*}
when $\tau_k$ is the first time the angel of power $c$ goes outside the $n$-dimensional sphere of radius $k$.
\end{thm}

\begin{thm}\label{prob}
Let $k$ be the positive integer and $\tau_k$ be the first time angel of power $c$ goes outside the $n$-dimensional sphere of radius $k$. Then
\begin{equation*}
    \mathbb{P}[\tau_k \le N] \le \frac{nc(c+1)N}{3k^2} .
\end{equation*}
\label{thm:hitting}
\end{thm}

\begin{cor}
If the hiding devil plays with the Gauss's circle strategy on a $2$-dimensional infinite grid graph, then,
\begin{equation*}
    \mathbb{P}[A^{2}_{k, c}] \ge 1 - \frac{4\pi c(c+1)}{3k}\left[ \sqrt{2}(c+1) + \frac{c(c+2)}{k} \right].
\end{equation*}
\label{cor:hitting}
\end{cor}


\section{Preliminaries}\label{sec3}

In the first subsection, we introduce a well-known problem called \emph{Gauss's circle}. We use the solution of this problem to count the number of turns that the devil needs to construct the hidden cage.

\subsection{Gauss's Circle}\label{gauss}
The Gauss's circle problem aims to count $N(r)$, the number of the vertices $(x, y) \in \mathbb{Z} \times \mathbb{Z}$ in a circle of radius $r>0$ centred at $(0, 0)$ (including the boundary). It can be approximated that
 \begin{equation*}
 	N(r) = \pi r^2 + O(r^{t})
 \end{equation*}
where $\frac{1}{2} \le t \le \frac{131}{208}$ with the lower bound was proved by Hardy 1915 \cite{hardy1915expression} while the upper bound was proved Huxley in 2000 \cite{huxley2002integer}. Interestingly, $N(r)$ can be expressed in many different ways. One of classical formulae was established by Hilbert and Cohn-Vossen \cite{hilbert2021geometry} as follows:
 \begin{equation*}
	 N(r)=1+4\sum _{{i=0}}^{\infty }\left(\left\lfloor {\frac  {r^{2}}{4i+1}}\right\rfloor -\left\lfloor {\frac  {r^{2}}{4i+3}}\right\rfloor \right).
 \end{equation*}
As the problem has been answered for the case of $2$ dimensional space, it is natural to ask further when we study on a general number of dimensions: 
\vskip 5 pt

\indent For a natural number $n \geq 2$ and a positive real number $r$, how many points (vertices) $x \in \mathbb{Z}^{n}$ are there that lie inside or on the surface of $L^{n}_{r}$, the $n$-dimensional sphere of radius $r$ centered at the origin? Namely, the generalized problem is to find $N_n(r)$, the number of points $x \in \mathbb{Z}^{n}$ such that $\|x\| \leq r$.
\vskip 5 pt

\indent However, the above question can be answered by approximating $N_n(r)$ with $Vol(S^{n}_{r})$, the volume of $S^{n}_{r}$, for which Smith and Vamanamurthy \cite{smith1989small} gave several proofs in their paper that:
\begin{equation*}
	Vol(S^{n}_{r})={\frac {\pi ^{\frac {n}{2}}}{\Gamma \left({\frac {n}{2}}+1\right)}}r^{n}.
\end{equation*}

\subsection{Jacobian for $n$-dimensional spherical coordinates}\label{Jac}
Blumension \cite{blumenson1960derivation} found the Jacobian when a multiple integral is transformed from rectangular to spherical coordinate systems of $n$-dimensions for an arbitrary $n \in \mathbb{N}$. 
In general, the equation for the sphere of radius $r$ in $n$-dimensional rectangular coordinate is
\begin{equation*}
    x_1^2 + x_2^2 + \dots + x_n^2 = r^2
\end{equation*}
where $x_i$'s are Cartesian coordinates.
The transformation for $n$-dimensional spherical coordinates is
\begin{align*}
x_1 &= r\cos{\phi_1}\\
x_2 &= r\sin{\phi_1}\cos{\phi_2}\\
x_3 &= r\sin{\phi_1}\sin{\phi_2}\cos{\phi_3}\\
&\quad \vdots\\
x_i &= r\sin{\phi_1}\sin{\phi_2}\dots\sin{\phi_{i-1}}\cos{\phi_i}\\
&\quad \vdots\\
x_{n-2} &= r\sin{\phi_1}\sin{\phi_2}\dots\sin{\phi_{n-3}}\cos{\phi_{n-2}}\\
x_{n-1} &= r\sin{\phi_1}\sin{\phi_2}\dots\sin{\phi_{n-2}}\cos{\theta}\\
x_{n} &= r\sin{\phi_1}\sin{\phi_2}\dots\sin{\phi_{n-2}}\sin{\theta}\\
\end{align*}
where $0 \le \phi_i \le \pi$ for $i \in \{1,2,3,\dots,n-2 \}$ and $0 \le \theta \le 2\pi$.\\
The Jacobian is a determinant of the $n$ by $n$ matrix of partial derivatives
\begin{equation*}
  J(r,\theta,\phi_1,\dots,\phi_{n-2}) = \begin{bmatrix}
\frac{\partial x_1}{\partial r} & \frac{\partial x_1}{\partial \theta} & \frac{\partial x_1}{\partial \phi_1} & \cdots & \frac{\partial x_1}{\partial \phi_{n-2}}\\
\frac{\partial x_2}{\partial r} & \frac{\partial x_2}{\partial \theta} & \frac{\partial x_2}{\partial \phi_1} & \cdots & \frac{\partial x_2}{\partial \phi_{n-2}}\\
\vdots & \vdots & \vdots & \ddots & \vdots \\
\frac{\partial x_n}{\partial r} & \frac{\partial x_n}{\partial \theta} & \frac{\partial x_n}{\partial \phi_1} & \cdots & \frac{\partial x_n}{\partial \phi_{n-2}}
\end{bmatrix} 
= r^{n-1}\prod_{k=1}^{n-2} \sin^{n-1-k}{\phi_k}.
\end{equation*}

\subsection{Central limit theorem}\label{CLT}


\begin{thm}[Multivariate central limit theorem]
    \label{MCLT}
    Let $\boldsymbol{\xi}_1,\boldsymbol{\xi}_2,\dots,\boldsymbol{\xi}_m$ be $m$ discrete independent and identically distributed $n$-dimensional random vectors with mean $\boldsymbol{\mu}$ and finite covariance matrix $\boldsymbol{\Sigma}$ and $\mathbf{X}_m$ be the sum of $\boldsymbol{\xi}_i$ for $i=1,2,\dots,m$. Then
    \begin{equation*}
        \frac{\mathbf{X}_m-m\boldsymbol{\mu}}{\sqrt{m}} \to \Phi_{\boldsymbol{\Sigma}}
    \end{equation*}
    where $\Phi_{\boldsymbol{\Sigma}}$ is an $\mathcal{N}_n(\mathbf{0},\boldsymbol{\Sigma})$ distribution.
\end{thm}

By using central limit theorem, the probability that the angel of power $c$ is caged inside the hollow sphere of radius $k$ after $N$ moves can be approximated by
\begin{equation*}
    P[A_{k,c}^{n}] = P[\|\mathbf{X}_{N}\|^{2} < k^2] \approx \int \cdots \int_{D} f_{Z}(x_1,x_2,\dots,x_n) dx_1 dx_2 \dots dx_n
\end{equation*}
where $D = \{(x_1,x_2,\dots,x_n) \in \mathbb{R}^n| x_{1}^{2} + x_{2}^{2} + \dots + x_{n}^{2} < k^2\}$ and $Z \sim \mathcal{N}_{n}(\mathbf{0},\boldsymbol{\Sigma})$. The covariance matrix, $\boldsymbol{\Sigma}$, can be calculated by
\begin{align*}
    \boldsymbol{\Sigma} &= \begin{bmatrix}
		Var(X_{N}^{(1)}) & 0 & \cdots & 0\\
		0 & Var(X_{N}^{(2)}) & \cdots & 0\\
		\vdots & \vdots & \ddots & \vdots\\
		0 & 0 & \cdots & Var(X_{N}^{(n)})
		\end{bmatrix}\\
		\\
		&= \begin{bmatrix}
		Var(\sum_{i=1}^{N} \xi_{i}^{(1)}) & 0 & \cdots & 0\\
		0 & Var(\sum_{i=1}^{N} \xi_{i}^{(2)}) & \cdots & 0\\
		\vdots & \vdots & \ddots & \vdots\\
		0 & 0 & \cdots & Var(\sum_{i=1}^{N} \xi_{i}^{(n)})
		\end{bmatrix}\\
		\\
		&= \begin{bmatrix}
		\frac{N(c^2+c)}{3} & 0 & \cdots & 0\\
		0 & \frac{N(c^2+c)}{3} & \cdots & 0\\
		\vdots & \vdots & \ddots & \vdots\\
		0 & 0 & \cdots & \frac{N(c^2+c)}{3}
		\end{bmatrix}.\\
\end{align*}

\section{Proof of Lemma \ref{dim}}\label{sec4}
For the equation
\begin{equation*}
	x_1 + x_2 + \dots + x_l = k
\end{equation*}
with $x_i \in \{-c,-c+1,\dots , 0, \dots, c-1, c\}$ for all $i = 1,2,\dots, l$, the number of its solutions is equal to the number of the solutions of the following equation
\begin{equation*}
	x_1 + x_2 + \dots + x_l = k+cl
\end{equation*}
with $x_i \in \{0, 1, \dots, 2c\}$ for all $i = 1,2,\dots, l$ which is equal to the coefficient of $x^{k+cl}$ in the expansion of $(1 + x + x^2 + x^3 + \dots + x^{2c})^l$. Thus, by binomial expansion, we get
\begin{align*}
	(1 + x + x^2 + &x^3 + \dots + x^{2c})^l = \frac{(1-x^{2c+1})^l}{(1-x)^l}\\
	&= \Bigg(\sum_{m=0}^{l}(-1)^m\binom{l}{m}(x^{2c+1})^m\Bigg)\Bigg(1+\binom{l}{1}x +\binom{l+1}{2}x^2 + \binom{l+2}{3}x^3 + \dots \Bigg).
\end{align*}
So the coefficient of $x^{k+cl}$ is
\begin{equation*}
	\sum_{m=0}^{l} (-1)^m \binom{l}{m}\binom{l-1+k+cl-m(2c+1)}{l-1}.
\end{equation*}

\section{Proof of Theorem \ref{submain1}}\label{sec5}
We distinguish two cases as detailed in Subsections \ref{dim1} and \ref{dim2}.

\subsection{1-dimensional infinite grid graph}\label{dim1}
If the angel of power $c$ is at $x \in \mathbb{Z}$, then, in the next turn, the angel can jump to $x-c, x-c+1, \dots, x + c - 1$ or $x+c$, each of which with the probability $\frac{1}{2c+1}$. So, the angel move has a discrete uniform distribution on $-c, -c+1, \dots, 0 , \dots, c-1,c$. In 1-dimensional infinite grid graph, the number of vertices that the devil has to destroy to cage angel is $N = 2c$. The devil can destroy the vertices distance $c^2$ away from the initial position of the angel $(X_0=0)$. The devil is guaranteed to capture the angel.
\noindent This proves the case when $n = 1$.

\begin{figure}[h]
\centering
\tikzset{every picture/.style={line width=0.75pt}} 

\begin{tikzpicture}[x=0.75pt,y=0.75pt,yscale=-1,xscale=1]

\draw    (53,119.99) -- (546.8,118.81) ;
\draw [shift={(549.8,118.8)}, rotate = 179.86] [fill={rgb, 255:red, 0; green, 0; blue, 0 }  ][line width=0.08]  [draw opacity=0] (8.93,-4.29) -- (0,0) -- (8.93,4.29) -- cycle    ;
\draw [shift={(50,120)}, rotate = 359.86] [fill={rgb, 255:red, 0; green, 0; blue, 0 }  ][line width=0.08]  [draw opacity=0] (8.93,-4.29) -- (0,0) -- (8.93,4.29) -- cycle    ;
\draw    (300.8,111.8) -- (300.8,127.8) ;
\draw  [color={rgb, 255:red, 0; green, 0; blue, 0 }  ,draw opacity=0.5 ][fill={rgb, 255:red, 0; green, 0; blue, 0 }  ,fill opacity=0.48 ] (101,110.8) -- (171,110.8) -- (171,129.8) -- (101,129.8) -- cycle ;
\draw  [color={rgb, 255:red, 0; green, 0; blue, 0 }  ,draw opacity=0.5 ][fill={rgb, 255:red, 0; green, 0; blue, 0 }  ,fill opacity=0.48 ] (430,109.8) -- (500,109.8) -- (500,128.8) -- (430,128.8) -- cycle ;
\draw   (100.8,135.8) .. controls (100.8,140.47) and (103.13,142.8) .. (107.8,142.8) -- (123.64,142.8) .. controls (130.31,142.8) and (133.64,145.13) .. (133.64,149.8) .. controls (133.64,145.13) and (136.97,142.8) .. (143.64,142.8)(140.64,142.8) -- (162.8,142.8) .. controls (167.47,142.8) and (169.8,140.47) .. (169.8,135.8) ;
\draw   (300.8,98) .. controls (300.8,93.33) and (298.47,91) .. (293.8,91) -- (246.3,91) .. controls (239.63,91) and (236.3,88.67) .. (236.3,84) .. controls (236.3,88.67) and (232.97,91) .. (226.3,91)(229.3,91) -- (178.8,91) .. controls (174.13,91) and (171.8,93.33) .. (171.8,98) ;
\draw   (429.8,97) .. controls (429.8,92.33) and (427.47,90) .. (422.8,90) -- (375.3,90) .. controls (368.63,90) and (365.3,87.67) .. (365.3,83) .. controls (365.3,87.67) and (361.97,90) .. (355.3,90)(358.3,90) -- (307.8,90) .. controls (303.13,90) and (300.8,92.33) .. (300.8,97) ;
\draw   (430.8,135.8) .. controls (430.8,140.47) and (433.13,142.8) .. (437.8,142.8) -- (453.64,142.8) .. controls (460.31,142.8) and (463.64,145.13) .. (463.64,149.8) .. controls (463.64,145.13) and (466.97,142.8) .. (473.64,142.8)(470.64,142.8) -- (492.8,142.8) .. controls (497.47,142.8) and (499.8,140.47) .. (499.8,135.8) ;

\draw (295,130) node [anchor=north west][inner sep=0.75pt]   [align=left] {$0$};
\draw (128,154) node [anchor=north west][inner sep=0.75pt]   [align=left] {$c$};
\draw (459,154) node [anchor=north west][inner sep=0.75pt]   [align=left] {$c$};
\draw (227,64) node [anchor=north west][inner sep=0.75pt]   [align=left] {$c^2$};
\draw (356,64) node [anchor=north west][inner sep=0.75pt]   [align=left] {$c^2$};

\end{tikzpicture}
\caption{Strategy for Devil in 1-dimensional infinite grid graph}
\end{figure}
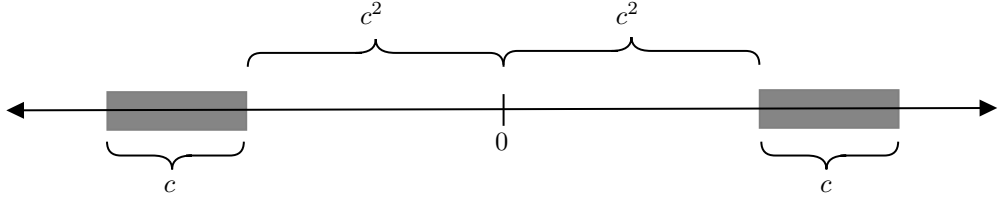

\subsection{2-dimensional infinite grid graph}\label{dim2}
If the angel of power $c$ is at $(x,y) \in \mathbb{Z}^2$, then, in the next turn, she can jump to any vertex $(x', y')$ which $\max\{|x - x'|, |y - y'|\} \leq c$, each of which with the probability $\frac{1}{(2c+1)^2}$. By Gauss's circle problem, the number of vertices that the devil has to destroy to cage the angel in $H^{2}_{k, c}$ is
\begin{align*}
	N = |H^{2}_{k, c}| &\approx \pi(k+c\sqrt{2})^2 - \pi(k)^2\\
	&= \pi(c\sqrt{2})(2k+c\sqrt{2}).
\end{align*}
We can see that the position of the angel can be calculated separately in $x$-axis and $y$-axis. Suppose $\boldsymbol{\xi}$ be the random vector for jumping in $x$-axis and $y$-axis.
\begin{align*}
	\boldsymbol{\xi}_i = \begin{bmatrix}
        \xi^{(1)}_{i} \\
        \xi^{(2)}_{i}
        \end{bmatrix}, \boldsymbol{\mu}_{\boldsymbol{\xi}_i} = \begin{bmatrix}
        0 \\
        0
        \end{bmatrix},
        \boldsymbol{\Sigma}_{\boldsymbol{\xi}_i} = \begin{bmatrix}
        \frac{c^2+c}{3} & 0\\
        0 & \frac{c^2+c}{3}
        \end{bmatrix}
\end{align*} 
The position of angel after $N$ moves starting at $(0,0)$ is the sum of $N$ identically independently distributed random variables:
\begin{align*}
	\mathbf{X}_N &= \begin{bmatrix}
        X^{(1)}_{N} \\
        X^{(2)}_{N}
        \end{bmatrix} = \begin{bmatrix}
        \sum_{i=0}^{N} \xi^{(1)}_{i} \\
        \sum_{i=0}^{N} \xi^{(2)}_{i}
        \end{bmatrix}\\
	\boldsymbol{\mu}_{\mathbf{X}_N} = \begin{bmatrix}
        0 \\
        0
        \end{bmatrix}, \quad \boldsymbol{\Sigma}_{\mathbf{X}_N} &= \begin{bmatrix}
        \frac{\sqrt{2} \pi c^2(c+1)(2k+c\sqrt{2})}{3} & 0\\
        0 & \frac{\sqrt{2} \pi c^2(c+1)(2k+c\sqrt{2})}{3}
        \end{bmatrix} 
\end{align*}
By Section \ref{CLT}, the probability that the angel of power $c$ is caged inside the hollow sphere of radius $k$ after $N$ moves can be approximated by
\begin{align*}
    P[A^{2}_{k, c}] &\approx \int \int_{D} f_{\mathbf{X}_N}(z) dx_1 dx_2\\
    &= \int \int_{D} \frac{exp(-\frac{1}{2}z^T\Sigma^{-1}z)}{2\pi\sqrt{|\Sigma|}} dx_1 dx_2\\
    &= \int \int_{D} \frac{exp(-\frac{(x^2+y^2)}{2\sigma^2})}{2\pi\sigma^2} dx_1 dx_2
\end{align*}
where $D = \{(x_1,x_2) \in \mathbb{R}^n| x_{1}^{2} + x_{2}^{2} < k^2\}$ and $\sigma^2 = \frac{\sqrt{2} \pi c^2(c+1)(2k+c\sqrt{2})}{3}$.
By applying double integral in polar coordinate, 
\begin{align*}
\int \int_{D} f_{\mathbf{X}_N}(z) dx_1 dx_2 = \int_{0}^{2\pi}\int_{0}^{k} f_{\mathbf{X}_N}(r) r drd\theta &= \int_{0}^{2\pi}\int_{0}^{k} \frac{exp(-\frac{(r^2)}{2\sigma^2})}{2\pi\sigma^2} r drd\theta\\
	&= \frac{1}{2\pi\sigma^2} \int_{0}^{2\pi}\int_{0}^{k} r \cdot exp(-\frac{(r^2)}{2\sigma^2})  drd\theta\\
	&= \frac{1}{2\pi\sigma^2}\Big(\int_{0}^{2\pi}1d\theta\Big) \Big(\int_{0}^{k} r \cdot exp(-\frac{(r^2)}{2\sigma^2})  dr \Big)\\
	&= 1 - exp\big(-\frac{k^2}{2\sigma^2}\big)
\end{align*}
\noindent Hence, for the value $k$ such that $P[A^{2}_{k, c}] > 1-\epsilon$, we have that $k$ satisfies the following inequality
\begin{equation}
	exp\big(-\frac{k^2}{2\sigma^2}\big) \le \epsilon.
	\label{2d-condition}
\end{equation}
\noindent By Subsections \ref{dim1} and \ref{dim2}, we prove Theorem \ref{submain1}.

\section{Proof of Theorems \ref{submain2} and \ref{submain3}}\label{sec6}
For $n \geq 3$, if the angel is at $x \in \mathbb{Z}^n$, then, in the next turn, the angel can jump to any vertex $y$ such that $d(x, y) \leq c$, each of which with the probability $\frac{1}{(2c+1)^n}$. By Gauss's circle problem, the number of vertices that the devil has to destroy to cage the angel in $H^{n}_{k, c}$ is
\begin{align*}
	N = |H^{n}_{k, c}| &\approx \frac{\pi^{\frac{n}{2}}[(k+c\sqrt{n})^n]}{\Gamma(\frac{n}{2}+1)} - \frac{\pi^{\frac{n}{2}}[k^n]}{\Gamma(\frac{n}{2}+1)}\\
	&= \frac{\pi^{\frac{n}{2}}[(k+c\sqrt{n})^n-k^n]}{\Gamma(\frac{n}{2}+1)}.
\end{align*}
We can see that the position of the angel can be calculated separately in each axis. Suppose $\boldsymbol{\xi}_i$ be the random vector for jumping of the angel.
\begin{align*}
	\boldsymbol{\xi}_i = \begin{bmatrix}
        \xi^{(1)}_{i} \\
        \xi^{(2)}_{i} \\
        \vdots \\
        \xi^{(n)}_{i}
        \end{bmatrix}, \boldsymbol{\mu}_{\boldsymbol{\xi}_i} = \begin{bmatrix}
        0 \\
        0 \\
        \vdots \\
        0
        \end{bmatrix},
        \boldsymbol{\Sigma}_{\boldsymbol{\xi}_i} = \begin{bmatrix}
		\frac{N(c^2+c)}{3} & 0 & \cdots & 0\\
		0 & \frac{N(c^2+c)}{3} & \cdots & 0\\
		\vdots & \vdots & \ddots & \vdots\\
		0 & 0 & \cdots & \frac{N(c^2+c)}{3}
        \end{bmatrix}
\end{align*} 
The position of the angel after $N$ moves starting at the origin is the sum of $N$ identically independently distributed random variables. Let $\mathbf{X}_N = (X^1_N,X^2_N,\dots,X^n_N)^T$. Thus,
\begin{align*}
	\mathbf{X}_N &= \begin{bmatrix}
        X^{(1)}_{N} \\
        X^{(2)}_{N} \\
        \vdots \\
        X^{(n)}_{N}
        \end{bmatrix} = \begin{bmatrix}
        \sum_{i=0}^{N} \xi^{(1)}_{i} \\
        \sum_{i=0}^{N} \xi^{(2)}_{i} \\
        \vdots \\
        \sum_{i=0}^{N} \xi^{(n)}_{i}
        \end{bmatrix}\\ \\
	\boldsymbol{\mu}_{\mathbf{X}_N} = \begin{bmatrix}
        0 \\
        0 \\
        \vdots \\
        0
        \end{bmatrix}, \quad \boldsymbol{\Sigma}_{\mathbf{X}_N} &=  \begin{bmatrix}
		\frac{N(c^2+c)}{3} & 0 & \cdots & 0\\
		0 & \frac{N(c^2+c)}{3} & \cdots & 0\\
		\vdots & \vdots & \ddots & \vdots\\
		0 & 0 & \cdots & \frac{N(c^2+c)}{3}
		\end{bmatrix}\\
\end{align*}
Let $\sigma^2 = \frac{c^2+c}{3}N$. By Section \ref{CLT}, the probability that the angel of power $c$ is caged inside the hollow sphere of radius $k$ after $N$ moves can be approximated by
\begin{align*}
    P[A^{n}_{k, c}] &\approx \int \dots \int_{D} f_{\mathbf{X}_N}(z) dx_1 dx_2 \dots dx_n\\
    &= \int \dots \int_{D} \frac{exp(-\frac{1}{2}z^T\Sigma^{-1}z)}{2\pi\sqrt{|\Sigma|}} dx_1 dx_2 \dots dx_n\\
    &= \int \dots \int_{D} \frac{exp(-\frac{(x_1^2+x_2^2+\dots+x_n^2)}{2\sigma^2})}{2\pi\sigma^2} dx_1 dx_2 \dots dx_n
\end{align*}
where $D = \{(x_1,x_2) \in \mathbb{R}^n| x_{1}^{2} + x_{2}^{2} + \dots + x_{n}^{2} < k^2\}$.
By the transformation for $n$-dimensional spherical coordinates (see Subsection \ref{Jac}), the probability that the angel of power $c$ will be in $H^{n}_{k, c}$ after $N$ moves starting from the origin is
\begin{align*}
	P[A^{n}_{k, c}] &\approx	\int_{0}^{2\pi}\int_{0}^{\pi}\cdots \int_{0}^{\pi}\int_{0}^{k} f_{Z}(r) r^{n-1}\prod_{k=1}^{n-2}\sin^{n-1-k}\phi_k drd\phi_1\cdots d\phi_{n-2}d\theta\\
	&= \int_{0}^{2\pi}\int_{0}^{\pi}\cdots \int_{0}^{\pi}\int_{0}^{k} \frac{exp(-\frac{r^2}{2\sigma^2})}{\sigma^n\sqrt{(2\pi)^n}} r^{n-1}\prod_{k=1}^{n-2}\sin^{n-1-k}\phi_k drd\phi_1\cdots d\phi_{n-2}d\theta\\
	&= \frac{1}{\sigma^n\sqrt{(2\pi)^n}}\Big(\int_{0}^{2\pi}1d\theta\Big)\Big(\int_{0}^{\pi}\sin^{n-2}\phi_1 d\phi_1\Big)\cdots \Big(\int_{0}^{\pi}\sin\phi_{n-2} d\phi_{n-2}\Big) \Big(\int_{0}^{k} r^{n-1} \cdot exp(-\frac{r^2}{2\sigma^2})  dr \Big).
\end{align*}
Then, by the technique from White \cite{white1960matrix}, we distinguish our integral $2$ cases according to the parity of $n$. 
When $n$ is even, we have that
\begin{align*}
	&\frac{1}{\sigma^n\sqrt{(2\pi)^n}}\Big(\int_{0}^{2\pi}1d\theta\Big)\Big(\int_{0}^{\pi}\sin^{n-2}\phi_1 d\phi_1\Big)\cdots \Big(\int_{0}^{\pi}\sin\phi_{n-2} d\phi_{n-2}\Big) \Big(\int_{0}^{k} r^{n-1} \cdot exp(-\frac{r^2}{2\sigma^2})  dr \Big)\\
	&= \frac{1}{\sigma^n\sqrt{(2\pi)^n}} (2 \pi) \Big(\frac{2^{\frac{n}{2}-1}\pi^{\frac{n}{2}-1}}{(n-2)!!}\Big)\Big(2^{\frac{n}{2}-1}\sigma^n(\frac{n}{2}-1)!-2^{\frac{n}{2}-1}\sigma^n(\frac{n}{2}-1)!exp(-\frac{k^2}{2\sigma^2})\sum_{d=0}^{\frac{n}{2}-1}\frac{k^{2d}}{2^d\sigma^{2d}d!}\Big)\\
	&= \frac{2^{\frac{n}{2}-1}({\frac{n}{2}-1})!}{(n-2)!!}\Big[1 - exp(-\frac{k^2}{2\sigma^2})\sum_{d=0}^{\frac{n}{2}-1}\frac{k^{2d}}{2^d\sigma^{2d}d!}\Big]\\
	&= 1 - exp(-\frac{k^2}{2\sigma^2})\sum_{d=0}^{\frac{n}{2}-1}\frac{k^{2d}}{2^d\sigma^{2d}d!}.
\end{align*}
\vskip 5 pt

\noindent Thus, when $n$ is even, we have
\begin{align*}
    P[A^{n}_{k, c}] \approx 1 - exp(-\frac{k^2}{2\sigma^2})\sum_{d=0}^{\frac{n}{2}-1}\frac{k^{2d}}{2^d\sigma^{2d}d!}.
\end{align*}

\noindent When $n$ is odd, we have that
\begin{align*}
	&\frac{1}{\sigma^n\sqrt{(2\pi)^n}}\Big(\int_{0}^{2\pi}1d\theta\Big)\Big(\int_{0}^{\pi}\sin^{n-2}\phi_1 d\phi_1\Big)\cdots \Big(\int_{0}^{\pi}\sin\phi_{n-2} d\phi_{n-2}\Big) \Big(\int_{0}^{k} r^{n-1} \cdot exp(-\frac{r^2}{2\sigma^2})  dr \Big)\\
	&= \frac{1}{\sigma^n\sqrt{(2\pi)^n}} (2 \pi)\Big(\frac{2^{\frac{n-1}{2}}\pi^{\frac{n-3}{2}}}{(n-2)!!}\Big)\Big((n-2)!! \sigma^n exp(-\frac{k^2}{2\sigma^2}) {\sum_{j=n}^{\infty}}{'}\frac{k^j}{\sigma^j j!!}\Big)\\
	&= \frac{\sqrt{2}}{\sqrt{\pi}}\Big[exp(-\frac{k^2}{2\sigma^2}) {\sum_{j=n}^{\infty}}{'}\frac{k^j}{\sigma^j j!!}\Big]
\end{align*}
where ${\sum_{j=n}^{\infty}}{'}$ is the summation from $j=n, n+2, n+4, \dots$
\vskip 5 pt

\noindent Thus, when $n$ is odd, we have
\begin{align*}
    P[A^{n}_{k, c}] \approx \frac{\sqrt{2}}{\sqrt{\pi}}\Big[exp(-\frac{k^2}{2\sigma^2}) {\sum_{j=n}^{\infty}}{'}\frac{k^j}{\sigma^j j!!}\Big].
\end{align*}

\subsection{$3$-dimensional infinite grid graph: proof of Theorem \ref{submain2}.}
In this case, we have 
$$N = \frac{4\pi}{3}\left[(k+c\sqrt{3})^3 - k^3 \right]$$ 
\noindent and 
$$\sigma = \sqrt{\frac{4\pi(c^2+c)}{9}\left[(k+c\sqrt{3})^3 - k^3 \right]}$$
\noindent Thus,
\begin{align*}
    \sigma &= \sqrt{\frac{4\pi(c^2+c)}{9}\left[(k+c\sqrt{3})^3 - k^3 \right]} \\
    &= \frac{2}{3}\sqrt{\pi c(c+1)\left[(3\sqrt{3}c)k^2+(9c^2)k+3\sqrt{3}c^3)\right]}\\
    &= \frac{2c}{3}\sqrt{\pi (c+1)\left[(3\sqrt{3})k^2+(9c)k+3\sqrt{3}c^2)\right]}.\\
\end{align*}
We get
\begin{equation*}
    \frac{2kc}{3}\sqrt{3\sqrt{3}\pi (c+1)} \le \sigma \le \frac{2kc}{3}\sqrt{3\sqrt{3}\pi (c+1)^3}
\end{equation*}

\noindent Further, 
\begin{align*}
	{\sum_{j=3}^{\infty}}{'}\frac{k^j}{\sigma^j j!!} 
	&= \sum_{l=1}^{\infty}\frac{l!2^l k^{(2l+1)}}{\sigma^{\left(2l+1\right)} (2l+1)!}\\
	&\le \sum_{l=1}^{\infty}\frac{l!2^l k^{(2l+1)}}{\left(\frac{2kc}{3}\sqrt{3\sqrt{3}\pi (c+1)}\right)^{2l+1} (2l+1)!}\\
        &= \sum_{l=1}^{\infty}\frac{l!3^{l+1} }{(2^{l+1})(c^{2l+1})\sqrt{3}\left(\sqrt{\sqrt{3}\pi (c+1)}\right)^{2l+1}(2l+1)!}\\	
        &= \sum_{l=1}^{\infty} \left(\frac{3}{2}\right)^{l+1}\left(\frac{1}{\sqrt{3}c^{2l+1}\left(\sqrt{\sqrt{3}\pi (c+1) }\right)^{2l+1}}\right)\left(\frac{l!}{(2l+1)!}\right)\\
        &= \frac{3}{2\sqrt{3}c\sqrt{\sqrt{3}\pi (c+1) }}\sum_{l=1}^{\infty}\left(\frac{3}{2\sqrt{3}\pi c^2(c+1)}\right)^l\left(\frac{l!}{(2l+1)!}\right)\\
        &\le \frac{1}{4\sqrt{3}c\sqrt{\sqrt{3}\pi (c+1) }}\left(\frac{3}{2\sqrt{3}\pi c^2(c+1) - 3}\right)\\
        &= \left(\frac{1}{16\sqrt{3}\pi c^2(c+1) }\right)^{\frac{1}{2}}\left(\frac{1}{2\pi c^2(c+1) - \sqrt{3}}\right)\\
\end{align*}
and
\begin{align*}
	{\sum_{j=3}^{\infty}}{'}\frac{k^j}{\sigma^j j!!} 
	&= \sum_{l=1}^{\infty}\frac{l!2^l k^{(2l+1)}}{\sigma^{\left(2l+1\right)} (2l+1)!}\\
	&\ge \sum_{l=1}^{\infty}\frac{l!2^l k^{(2l+1)}}{\left(\frac{2kc}{3}\sqrt{3\sqrt{3}\pi (c+1)^3}\right)^{2l+1} (2l+1)!}\\
        &= \sum_{l=1}^{\infty}\frac{l!3^{l+1} }{(2^{l+1})(c^{2l+1})\sqrt{3}\left(\sqrt{\sqrt{3}\pi (c+1)^3}\right)^{2l+1} (2l+1)!}\\	
        &\ge \frac{3^2}{2^2 \sqrt{3} c^3 \left(\sqrt{\sqrt{3}\pi (c+1)^3}\right)^{3} 3!}\\
        &= \left(\frac{1}{2c(c+1)}\right)^3\left(\frac{\sqrt{3}}{3\pi^3(c+1)^3}\right)^{\frac{1}{2}}.
\end{align*}
We can see that ${\sum_{j=3}^{\infty}}{'}\frac{k^j}{\sigma^j j!!}$ doesn't depend on $k$ and bounded by
\begin{equation*}
    \left(\frac{1}{2c(c+1)}\right)^3\left(\frac{\sqrt{3}}{3\pi^3(c+1)^3}\right)^{\frac{1}{2}} \le  {\sum_{j=3}^{\infty}}{'}\frac{k^j}{\sigma^j j!!} \le \left(\frac{1}{16\sqrt{3}\pi c^2(c+1) }\right)^{\frac{1}{2}}\left(\frac{1}{2\pi c^2(c+1) - \sqrt{3}}\right).
\end{equation*}
\vskip 5 pt

\noindent Thus, $P[A^{3}_{k, c}]$ is bounded above by
\begin{equation*}
    \frac{\sqrt{2}}{\sqrt{\pi}}exp\left(-\frac{k^2}{2\sigma^2}\right){\sum_{j=3}^{\infty}}{'}\frac{k^j}{\sigma^j j!!} \le exp\left( -\frac{\sqrt{3}}{8\pi c^2(c+1)^3}\right)\left(\frac{1}{8\sqrt{3}\pi^2 c^2(c+1) }\right)^{\frac{1}{2}}\left(\frac{1}{2\pi c^2(c+1) - \sqrt{3}}\right)
\end{equation*}
and is bounded below by
\begin{equation*}
   \frac{\sqrt{2}}{\sqrt{\pi}}exp\left(-\frac{k^2}{2\sigma^2}\right){\sum_{j=3}^{\infty}}{'}\frac{k^j}{\sigma^j j!!} \ge exp\left(-\frac{\sqrt{3}}{8\pi c^2(c+1)}\right)\left(\frac{1}{2c(c+1)}\right)^3\left(\frac{2\sqrt{3}}{3\pi^4(c+1)^3}\right)^{\frac{1}{2}}.
\end{equation*}
Since $S_l = {\sum_{j=3}^{l}}{'}\frac{k^j}{\sigma^j j!!}$ is an increasing sequence., it follows that $\frac{\sqrt{2}}{\sqrt{\pi}}exp\left(-\frac{k^2}{2\sigma}\right){\sum_{j=3}^{\infty}}{'}\frac{k^j}{\sigma^j j!!}$ converges to some non-zero constant when $c$ is fixed. This proves Theorem \ref{submain2}.

\subsection{$n$-dimensional infinite grid graphs when $n \geq 4$: proof of Theorem \ref{submain3}.}
It can be observed by Gauss's circle problem that, in an $n$-dimensional infinite grid graph, the number of vertices that the devil has to destroy is $O(k^{n-1})$. So, $\sigma^2 = O(k^{n-1})$.
\vskip 5 pt

\indent We first assume that $n$ is even and $n \ge 4$. 
Thus,
\begin{align*}
	1-exp(-\frac{k^2}{2\sigma^2})\sum_{d=0}^{\frac{n}{2}-1}\frac{k^{2d}}{2^d\sigma^{2d}d!} &= 1- exp(-\frac{O(k^2)}{2O(k^{n-1})})\sum_{d=0}^{\frac{n}{2}-1}\frac{O(k^{2d})}{2^dO(k^{nd-d})d!}\\
	&= 1- exp(-O(k^{-n+3}))\sum_{d=0}^{\frac{n}{2}-1}O(k^{-nd+3d})\\
	&= 1 - 1\\
	&= 0
\end{align*}

\indent We finally assume that $n$ is odd and $n \ge 5$. Thus,
\begin{align*}
	\frac{\sqrt{2}}{\sqrt{\pi}}\Big[exp(-\frac{k^2}{2\sigma^2}) {\sum_{j=n}^{\infty}}{'}\frac{k^j}{\sigma^j j!!}\Big] &= \frac{\sqrt{2}}{\sqrt{\pi}}\Big[exp(-O(k^{-n+3})){\sum_{j=n}^{\infty}}{'}\frac{O(k^j)}{O(k^{j\frac{n-1}{2}})}\Big]\\
	&= \frac{\sqrt{2}}{\sqrt{\pi}}\Big[O(k^{-\frac{n^2-3n}{2}})\Big]\\
	&= O(k^{-\frac{n^2-3n}{2}})\\
	&= o(1)
\end{align*}
when $k \to\infty$. This proves Theorem \ref{submain3}.

\section{Proofs of Theorems \ref{stopping}, \ref{thm:hitting} and Corollary \ref{cor:hitting}}\label{sec7}

\subsection{Proof of Theorem \ref{stopping}}
Let $X_n$ be the position of the angel of power $c$. $X_n$ is the Markov chain on the two dimensional integer lattice with the following transition probabilities:
\begin{align*}
    \mathbb{P}\left( \mathbf{X}_{N+1} = (X^{(1)}_N+\xi^{(1)}_{N+1},X^{(2)}_N+\xi^{(2)}_{N+1},\dots,X^{(n)}_N+\xi^{(n)}_{N+1})^T | \mathbf{X}_{N} = (X^{(1)}_N,X^{(2)}_N,\dots,X^{(n)}_N)^T \right) = \frac{1}{(2c+1)^n}
\end{align*}
for all $\xi^{(m)}_{N+1} \in \{-c, -c+1, \dots, 0 \dots, c-1, c\}$.\\
Let $\boldsymbol{\xi}_i$ be i.i.d. random $n$-dimensional vectors with the distribution
\begin{align*}
    \mathbb{P}\left( \boldsymbol{\xi}_i = (\xi^{(1)}_{i},\xi^{(2)}_{i},\dots,\xi^{(n)}_{i})^T  \right) = \frac{1}{(2c+1)^n}
\end{align*}
for all $\xi^{(m)}_{i} \in \{-c, -c+1, \dots, 0 \dots, c-1, c\}$. We can see that 
\begin{equation*}
    \mathbf{X}_0 = \mathbf{0} \textnormal{ and } \mathbf{X}_{N} = \mathbf{X}_{N-1} + \boldsymbol{\xi}_N = \sum_{k=1}^{N} \boldsymbol{\xi}_i.
\end{equation*}
Note that $\mathbb{E}(\boldsymbol{\xi}_i) = \mathbf{0}$. Next we want to find the variance of $\boldsymbol{\xi}_i$,
\begin{align*}
    V&ar(\boldsymbol{\xi}_i)\\ &= \mathbb{E}(|\boldsymbol{\xi}_i|^2) \\&= \sum_{\boldsymbol{\xi}_i \in \Omega} \frac{1}{(2c+1)^n}|\boldsymbol{\xi}_i|^2 \\
    &= \sum_{\xi^{(1)}_{i} = -c}^{c} \sum_{\xi^{(2)}_{i} = -c}^{c} \cdots \sum_{\xi^{(n)}_{i} = -c}^{c} \frac{1}{(2c+1)^n}\left( (\xi^{(1)}_{i})^2 + (\xi^{(2)}_{i})^2 + \dots + (\xi^{(n)}_{i})^2 \right)\\
    &= \frac{1}{(2c+1)^n} \sum_{\xi^{(1)}_{i} = -c}^{c} \sum_{\xi^{(2)}_{i} = -c}^{c} \cdots \sum_{\xi^{(n-1)}_{i} = -c}^{c} \left( \sum_{\xi^{(n)}_{i} = -c}^{c}(\xi^{(1)}_{i})^2 + \sum_{\xi^{(n)}_{i} = -c}^{c}(\xi^{(2)}_{i})^2 + \dots + \sum_{\xi^{(n)}_{i} = -c}^{c}(\xi^{(n)}_{i})^2 \right)\\
    &= \frac{1}{(2c+1)^n} \sum_{\xi^{(1)}_{i} = -c}^{c} \sum_{\xi^{(2)}_{i} = -c}^{c} \cdots \sum_{\xi^{(n-1)}_{i} = -c}^{c} \left( (2c+1)(\xi^{(1)}_{i})^2 + (2c+1)(\xi^{(2)}_{i})^2 + \dots + \frac{c(c+1)(2c+1)}{3} \right)\\
    &\hspace{5pt} \vdots\\
    &= \frac{1}{(2c+1)^n} \left(\frac{c(c+1)(2c+1)^n}{3} + \frac{c(c+1)(2c+1)^n}{3} + \dots + \frac{c(c+1)(2c+1)^n}{3}\right) \\
    &= \frac{nc(c+1)}{3}
\end{align*}
Now, we show that $\mathbf{M}_N = |\mathbf{X}_N|^2 - \frac{nc(c+1)}{3} N$ is martingale.
\begin{align*}
    \mathbb{E}\left(\mathbf{M}_{N+1} | \mathcal{F}_N \right) &= \mathbb{E}\left(|\mathbf{X}_{N+1}|^2 - \frac{nc(c+1)}{3} (N+1) | \mathcal{F}_N \right)\\
    &= \mathbb{E}\left(|\mathbf{X}_{N}+\boldsymbol{\xi}_{N+1}|^2 - \frac{nc(c+1)}{3} (N+1) | \mathcal{F}_N \right)\\
    &= \mathbb{E}\left(|\mathbf{X}_{N}|^2 + 2\mathbf{X}_N\cdot \boldsymbol{\xi}_{N+1} + |\boldsymbol{\xi}_{N+1}|^2  - \frac{nc(c+1)}{3}N - \frac{nc(c+1)}{3} | \mathcal{F}_N \right)\\
    &= |\mathbf{X}_{N}|^2 - \frac{nc(c+1)}{3}N + \mathbb{E}\left(|\boldsymbol{\xi}_{N+1}|^2| \mathcal{F}_N \right) - \frac{nc(c+1)}{3}\\
    &= |\mathbf{X}_{N}|^2 - \frac{nc(c+1)}{3}N \\
    &= \mathbf{M}_N
\end{align*}
By Optional Stopping Theorem,
\begin{align*}
    0 = \mathbb{E}(|\mathbf{X}_0|^2) =\mathbb{E}(\mathbf{M}_0) &= \mathbb{E}(\mathbf{M}_{\tau_k})\\
    &= \mathbb{E}\left(|\mathbf{X}_{\tau_k}|^2 - \frac{nc(c+1)}{3}\tau_k\right)\\
    &= \mathbb{E}\left(|\mathbf{X}_{\tau_k}|^2\right) - \frac{nc(c+1)}{3}\mathbb{E}\left(\tau_k\right)
\end{align*}
Note that
\begin{equation*}
    k^2 \le |\mathbf{X}_{\tau_k}|^2 \le (k+c\sqrt{n})^2
\end{equation*}
We get
\begin{equation*}
    \frac{3}{nc(c+1)}k^2 \le  \mathbb{E}(\tau_k) \le \frac{3}{nc(c+1)}(k+c\sqrt{n})^2.
\end{equation*}

\subsection{Proof of Theorem \ref{thm:hitting}}
Let $\boldsymbol{\xi}_i$ be i.i.d. random $n$-dimensional vectors with the distribution
\begin{align*}
    \mathbb{P}\left( \boldsymbol{\xi}_i = (\xi^{(1)}_{i},\xi^{(2)}_{i},\dots,\xi^{(n)}_{i})^T  \right) = \frac{1}{(2c+1)^n}
\end{align*}
for all $\xi^{(m)}_{i} \in \{-c, -c+1, \dots, 0 \dots, c-1, c\}$. Define $ \mathbf{X}_{N} = \sum_{i=1}^{N} \boldsymbol{\xi}_i$ to be the position of the Angel of power $c$ in turn $N$. Then, by Kolmogorov's inequality,
\begin{equation*}
    \mathbb{P}[\tau_k \le N] = \mathbb{P}[\max_{1 \le i \le N} | \mathbf{X}_i| \ge k] \le \frac{Var[ \mathbf{X}_N]}{k^2}
\end{equation*}
The variance $Var[ \mathbf{X}_N] = \sum_{i=1}^{N}Var[\boldsymbol{\xi}_i] = \frac{nc(c+1)N}{3}$.
This complete the proof.

\subsection{Proof of Corollary \ref{cor:hitting}}
$H^{2}_{k, c}$ is defined to be $L^{2}_{k + \sqrt{2}c} - L^{2}_{k}$ where $L^{2}_{k}$ is the number of vertices $(x,y) \in \mathbb{Z}\times \mathbb{Z}$ in a circle of radius $r$ centered at $(0,0)$. We can get the upper bound of the $H^{2}_{k, c}$ by
\begin{align*}
    |H^{2}_{k, c}| &= N(k + \sqrt{2}c) - N(k)\\
    &= [\pi(k + \sqrt{2}c)^2 + E(k + \sqrt{2}c)] - [\pi(k)^2 + E(k)]\\
    &\le \pi(k + \sqrt{2}c)^2 + |E(k + \sqrt{2}c)| - \pi k^2\\
    &= \pi (\sqrt{2}c)(2k + \sqrt{2}c) + 2\sqrt{2}\pi(k + \sqrt{2}c)\\
    &= 2\sqrt{2}\pi k(c+1) + 2\pi c(c+2)
\end{align*}
where $E(k)$ is the error.

The probability that the Devil can entrap the drunk Angel is $\mathbb{P}[A^{2}_{k, c}]$. This probability is the probability when the Devil plays optimal. Therefore,
\begin{align*}
    \mathbb{P}[A^{2}_{k, c}] &\ge  \mathbb{P}[\tau_k > |H^{2}_{k, c}|]\\
    &= 1-  \mathbb{P}[\tau_k \le |H^{2}_{k, c}|]\\
    &\ge 1- \frac{2c(c+1)|H^{2}_{k, c}|}{3k^2}\\
    &\ge 1- \frac{2c(c+1)(2\sqrt{2}\pi k(c+1) + 2\pi c(c+2))}{3k^2}\\
    &= 1- \frac{4\sqrt{2}\pi c (c+1)^2}{3k} + \frac{4\pi c^2 (c+1)(c+2)}{3k^2}\\
    &= 1- \frac{4\pi c(c+1)}{3k}\left[ \sqrt{2}(c+1) + \frac{c(c+2)}{k} \right]
\end{align*}

\section{Numerical Results}\label{numerical}
In this section, we present the \textit{Monte Carlo simulations} of the angel move. All the codes are uploaded in https://github.com/nuttanon19701/DAnHD.git. Recall that $N$ is the number of turns that the angel can move which is the number of vertices that the devil has to destroy to cage the angel in $H^{n}_{k, c}$. 
\vskip 5 pt

\indent Table \ref{tab:1} shows the results when the game is played on a $2$-dimensional infinite grid graph. We performed 100,000 simulations for each $\epsilon \in \{0.5, 0.1, 0.01\}$ and $c \in \{1, 3, 10\}$. The inner radius $k$ of $H^{2}_{k, c}$ that the devil tries to entrap the angel can be bounded below by \eqref{2d-condition}.

\begin{table}[h]
\centering
\begin{tabular}{|l|ccc|ccc|ccc|}
\hline
$\epsilon$        & \multicolumn{3}{c|}{0.5}                                             & \multicolumn{3}{c|}{0.1}                                              & \multicolumn{3}{c|}{0.01}                                             \\ \hline
$c$              & \multicolumn{1}{c|}{1}      & \multicolumn{1}{c|}{3}       & 10      & \multicolumn{1}{c|}{1}       & \multicolumn{1}{c|}{3}       & 10      & \multicolumn{1}{c|}{1}       & \multicolumn{1}{c|}{3}       & 10      \\ \hline
$k$              & \multicolumn{1}{c|}{10}     & \multicolumn{1}{c|}{151}     & 4525    & \multicolumn{1}{c|}{29}      & \multicolumn{1}{c|}{495}     & 15013   & \multicolumn{1}{c|}{57}      & \multicolumn{1}{c|}{986}     & 30017   \\ \hline
success   rate & \multicolumn{1}{c|}{0.53128} & \multicolumn{1}{c|}{0.50425} & 0.50104 & \multicolumn{1}{c|}{0.90603} & \multicolumn{1}{c|}{0.90064} & 0.90045 & \multicolumn{1}{c|}{0.99133} & \multicolumn{1}{c|}{0.99028} & 0.99017 \\ \hline
\end{tabular}
\caption{Numerical result for $2$-dimensional infinite grid graph}
\label{tab:1}
\end{table}

\begin{figure}[h]
\centering
\subcaptionbox{$\epsilon = 0.5$}
{\includegraphics[width=0.33\textwidth]{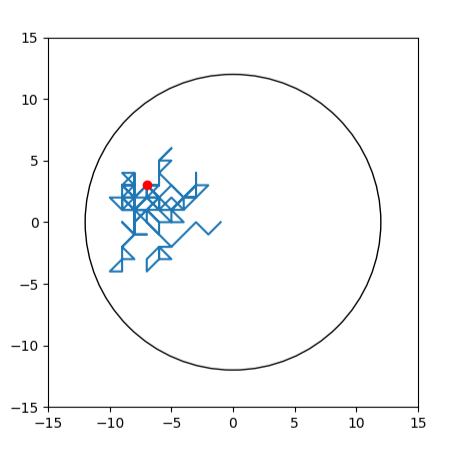}}%
\hfill
\subcaptionbox{$\epsilon = 0.1$}
{\includegraphics[width=0.33\textwidth]{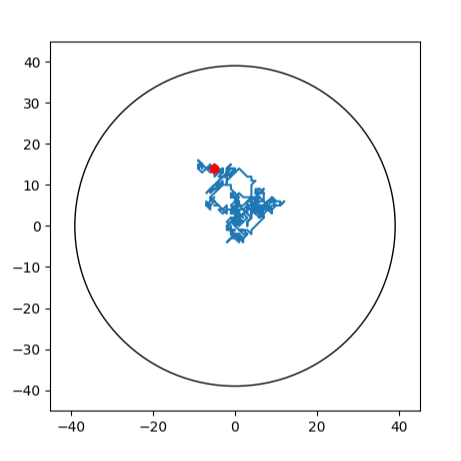}}%
\hfill
\subcaptionbox{$\epsilon = 0.01$}
{\includegraphics[width=0.33\textwidth]{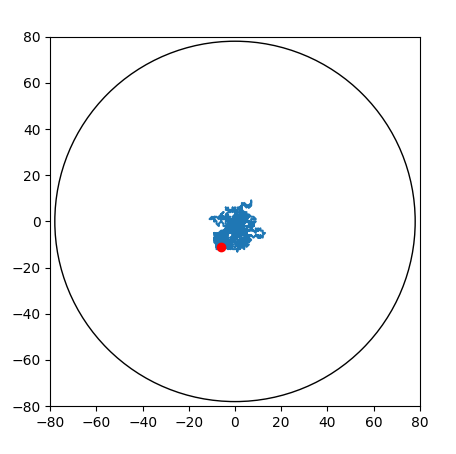}}%
\caption{The angel's footprints}
\label{fig:path}
\end{figure}

\noindent The results in Table \ref{tab:1} agree with Theorem \ref{submain1} that for any $c \in \mathbb{N}$ and $\epsilon > 0$, there exists $k$ such that $P[A^{2}_{k, c}] > 1 - \epsilon$. We further plot examples of the angel footprints in Figure \ref{fig:path} when $\epsilon \in \{0.5, 0.1, 0.01\}$ and $c = 1$. In the figure, the angel starts at the center of the circle and stops at the red vertex. It can be observed from Figure \ref{fig:path} that the angel travels closer to the origin than the boundary of $H^{2}_{k, 1}$ when $k$ is larger, when $\epsilon $ is smaller.
\vskip 10 pt

\vskip 15 pt

\begin{table}[h]
\centering
\begin{tabular}{|l|cccccc|}
\hline
dimension      & \multicolumn{3}{c|}{3}                                                                & \multicolumn{3}{c|}{4}                                       \\ \hline
$c$              & \multicolumn{6}{c|}{1}                                                                                                                               \\ \hline
$k$              & \multicolumn{1}{c|}{10}      & \multicolumn{1}{c|}{50} & \multicolumn{1}{c|}{100}     & \multicolumn{1}{c|}{5}        & \multicolumn{1}{c|}{10} & 25 \\ \hline
success   rate & \multicolumn{1}{c|}{0.00390} & \multicolumn{1}{c|}{0.00485}   & \multicolumn{1}{c|}{0.00395} & \multicolumn{1}{c|}{0.00001} & \multicolumn{1}{c|}{0.00000}  & 0.00000  \\ \hline
\end{tabular}
\caption{Numerical results for $n$-dimensional infinite grid graphs when $n \in \{3, 4\}$}
\label{tab:2}
\end{table}
\vskip 15 pt

\indent When the number of dimension is $3$, the results in Table \ref{tab:2} agree with Theorem \ref{submain2} that $P[A^{3}_{k, c}]$ is low but not $0$ even when $k$ is large. When the number of dimension is $4$, the results in Table \ref{tab:2} agree with Theorem \ref{submain3} that $P[A^{4}_{k, c}] \to 0$ when $k$ is large.
\vskip 5 pt

\newpage

\vspace*{\fill}
\begingroup
\centering

\begin{figure}[h]
\centering
\subcaptionbox{$2$-dimensional infinite grid graph}
{\includegraphics[width=0.50\textwidth]{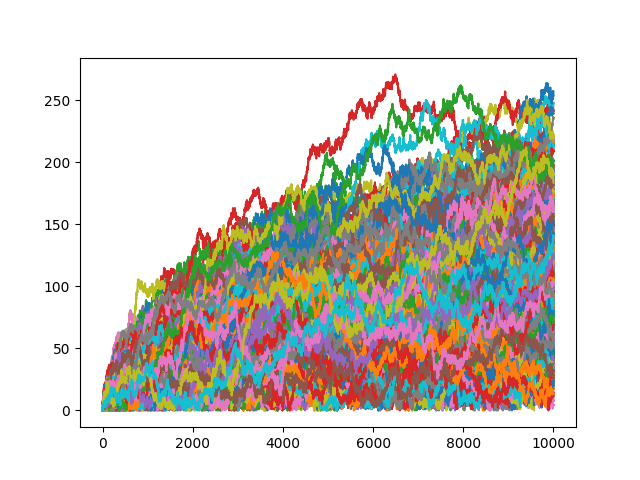}}%
\hfill
\subcaptionbox{$3$-dimensional infinite grid graph}
{\includegraphics[width=0.50\textwidth]{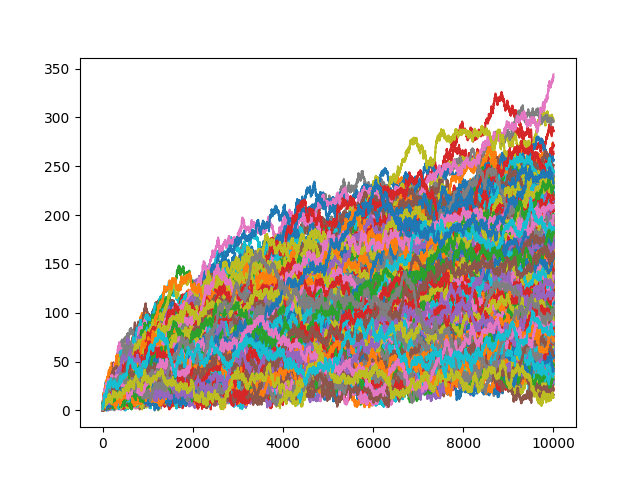}}%
\hfill
\subcaptionbox{$4$-dimensional infinite grid graph}
{\includegraphics[width=0.50\textwidth]{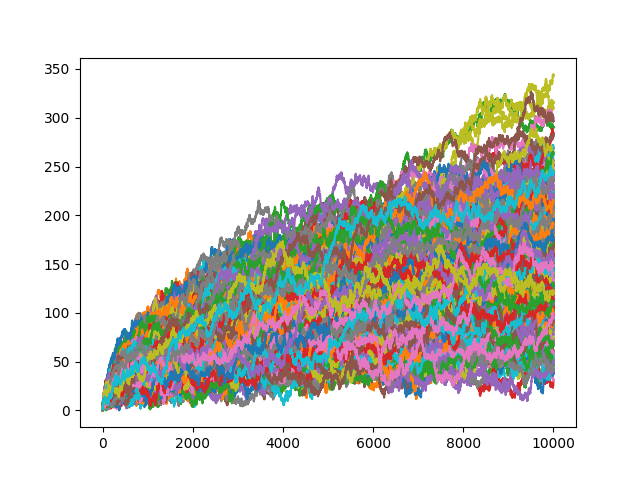}}%
\hfill
\subcaptionbox{$5$-dimensional infinite grid graph}
{\includegraphics[width=0.50\textwidth]{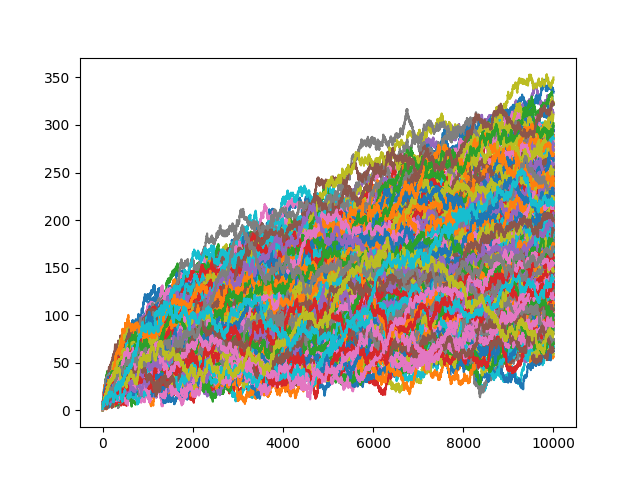}}%
\caption{Average distances of the angel between $0$ to $10,000$ steps}
\label{fig:monte carlo}
\end{figure}
\vskip 15 pt

\endgroup
\vspace*{\fill}

\indent Further, Figure \ref{fig:monte carlo} presents average distances ($y$-axis) of the angel from the origin between $0$ to $10,000$ moves ($x$-axis). The results show that the moving trend is farther from the origin when the number of dimension increases.

\newpage
\begin{figure}[h]
\centering
\subcaptionbox{$2$-dimensional infinite grid graph}
{\includegraphics[width=0.50\textwidth]{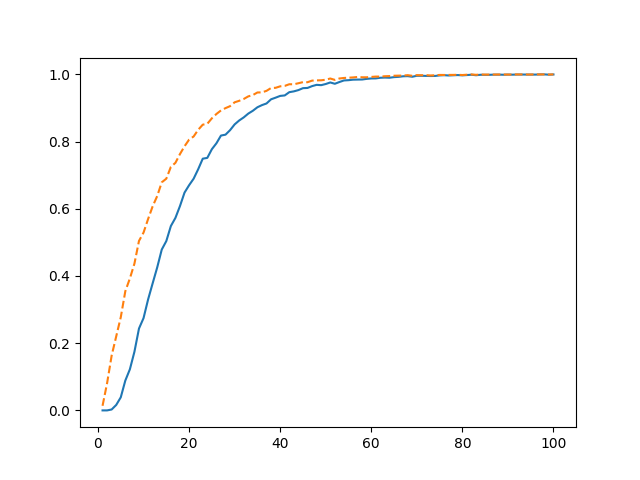}}%
\hfill
\subcaptionbox{$3$-dimensional infinite grid graph}
{\includegraphics[width=0.50\textwidth]{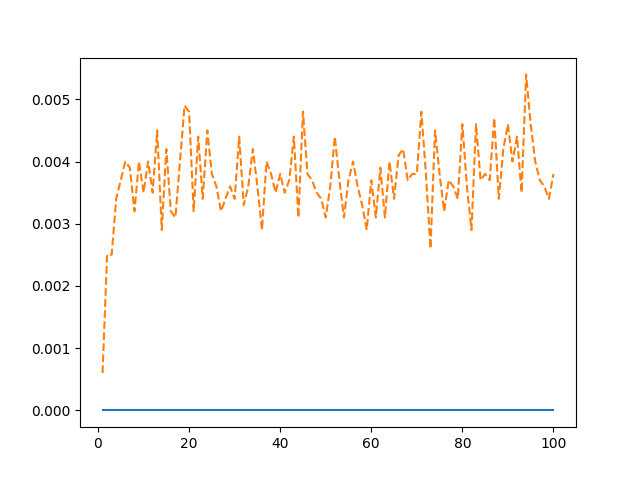}}%

\caption{The probabilities that the angel is in $H^{n}_{c, k}$ when $n \in \{2, 3\}$ and $k \in \{1, ..., 100\}$}
\label{fig:monte carlo1}
\end{figure}

\indent Finally, Figure \ref{fig:monte carlo1} presents Monte Carlo simulation of the probability ($y$-axis) that angel of power $c= 1$ is inside $H^{n}_{k, c}$ after $N$ moves (dashed) and the probability ($y$-axis) that the angel has not been outside $H^{n}_{k, c}$ (solid). The simulation is over 1,000,000 times with 10,000 of each $k \in \{1, ..., 100\}$ ($x$-axis). 

\bibliographystyle{abbrv} 
\bibliography{refs} 

\end{document}